# SEVERAL TYPICAL PARADIGMS OF INDUSTRIAL BIG DATA APPLICATION

Hu Shaolin, Zhang Qinghua, Su Naiquan and Li Xiwu

Guangdong University of Petrochemical Technology,
Maoming, Guangdong, China

*ABSTRACT*

*Industrial big data is an important part of big data family, which has important application value for industrial production scheduling, risk perception, state identification, safety monitoring and quality control, etc. Due to the particularity of the industrial field, some concepts in the existing big data research field are unable to reflect accurately the characteristics of industrial big data, such as what is industrial big data, how to measure industrial big data, how to apply industrial big data, and so on. In order to overcome the limitation that the existing definition of big data is not suitable for industrial big data, this paper intuitively proposes the concept of big data cloud and the 3M (Multi-source, Multi-dimension, Multi-span in time) definition of cloud-based big data. Based on big data cloud and 3M definition, three typical paradigms of industrial big data applications are built, including the fusion calculation paradigm, the model correction paradigm and the information compensation paradigm. These results are helpful for grasping systematically the methods and approaches of industrial big data applications.*

## 1. INTRODUCTION

Since the concept of big data was accepted by the Chinese in 2013, a big wave of big data has emerged across the country. On the one hand, the government has raised big data into a national development strategy and established a series of big data research institutes or big data centers[1]. On the other hand, for ordinary people, big data appear in communication in the way of common language in daily life. As mentioned in [1], the field of big data plays a vital role in various fields. But, in many different fields, the big data is just a term for massive data sets having large amounts of data, more varied and complex structure with the difficulties of analysing, storing and visualizing for further processes or results[1].

Although there are still many fundamental problems to be solved, big data technologies and applications are still in progress. In fact, the Big data approach and its applications are very extensive in various fields. Paper [1] presented an explanation of these applications such as telecommunication, business process management and human resources management. Paper [2] introduced and analysed the application trend of big data in power industry. In this paper, the concept definition and several typical use paradigms of big data in petrochemical industry and other industries will be discussed. It is worth mentioning that these basic problems, such as what is big data and how to process practical big data, affect the enthusiasm of technical personnel to innovate big data technology.

Over the years, although the concept of big data has been receiving increasing attention, some basic issues have not been solved until now[3,4], for example, what is big data, how to measure the big data, how to analyse and process the big data, and how to apply the big data. Especially in the industrial field, such as petrochemical industry, although the term "big data" has been widely mentioned, but there is no specific details that can be easily accepted by technical personnel. It is even more difficult for technical personnel to grasp the occasions and how to make full use of these "big data".

Up to now, the big data definition widely adopted by almost everyone is the descriptive definition of "4V" (Volume-large quantified volume, Velocity-fast velocity, Variety-diversified variety, Value-low value density, short as 4V). However, parts of 4V are inappropriate or ambiguous[5-6]. For one example, in a long-running actual petrochemical system, a single-dimensional long-term sequence sampling data collected by a sensor for a long time, although it may be very large, is not enough to constitute the so-called "big data" in the field of big data research. For the other example, low value density is a concept that is very difficult to understand. Maybe some big data does not have much value for big data. The key points depend on which perspective you look at, but it's not that big data necessarily has a low density of specific data.

In order to overcome the limitations of 4V definition of big data, and to more fully and accurately describe the characteristics of big data, the second section of this paper will present a new set of descriptive definitions of big data from a new perspective, and explore the measurement of big data volume. Based on the new definition of big data, the third section of this paper will briefly describe several typical paradigms of big data applications. These typical paradigms will help us understand the actual use of big data and how to implement these practical applications. At the end of this paper, several conclusions will be refined.

## 2. MEASURABLE DEFINITION FOR INDUSTRIAL BIG DATA

As we all know, "big" is a vague concept without clear boundaries. In other words, the "Big" is a vague qualitative descriptive adjective. What is "big" data and how to measure whether the "big" data is big or not big enough? The "4V" definition for big data does not rigorously tell us what the big data is, and how to determine whether the data is big data or not. The fuzziness of big data concept directly affects the exploration and application of big data by engineering and technical personnel.

In order to overcome some limitations about the 4V definition for big data[2,4], this section presents a new measurable definition that can be used to measure a data set is big or not.

Generally, the manufacturing or production process may continue for a period of time, the production process may be repeated again and again, and the products may be affected by various internal and external factors such as raw materials, environment, and processing disturbances, etc. So, the data sources of big data are extensive and have the certain time spans. In other words, the big data is a data cluster made up of many different data as well as data set:

$$S(\omega \in \Omega, t \in T, n \in N) \quad (1)$$

where the set $\Omega$ is the collection set of various sources, the set $T$ is the collection set of various time passages, the set $N$ is the collection set of various data dimensions.

In recent years, the concept of cloud has been widely accepted. Industrial big data is like a data cloud floating over factories. The cluster of big data described above can be represented graphically, as shown in Fig. 1.

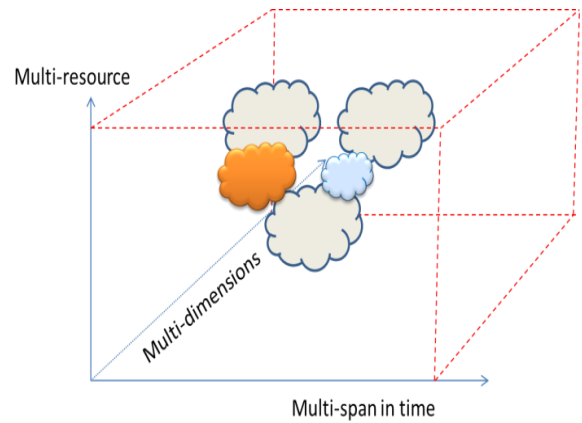

Figure 1. Big Data Cluster Diagram

Figure 1 shows that the big data cluster is very like several clouds floating in the air. Some of the clouds can be separated, and some clouds may be partially overlapped with other clouds. And, it can be seen from Figure 1 that big data has "3M" features: multi-source, multi-dimension and multi-span in time.

(1) Multiple sources. In the field of industrial production and processing, data can come from different scenarios, different production periods, different production links, different objects, and different incentives. For example, in the petrochemical production process, data can come from the working conditions of petrochemical equipment, the production process of petrochemical products, the external environment of the petrochemical production process, the internal state of petroleum refining, the quality of petroleum raw materials, as well as personal factors, etc. All these data are important for the safety management of petrochemical processes. All these data come from different sources, just like a colourful cloud floating in the air, and they come together to form big data of petrochemical industry. In the petrochemical production management process, each cloud has its own rationality and necessity.

(2) Multi-dimension: When we observe the real thing, we usually examine it from different perspectives and levels. Similarly, a petrochemical production process is taken as an example. In the production process, a large number of sensor devices are set up. Each sensor is like a person's eyes, observing a certain side of the chemical process, and obtaining information which display the states of the process partially. Observation information obtained from all different sides is brought together to form a set of high-dimensional data (even ultra-high-dimensional data). Although single-dimensional data is not enough to describe the overall change, the multi- dimensional data flow can be used to completely describe the changing process of petrochemical production.

(3) Multi-span in time. Industrial production is an ongoing process. In a sense, actions performed at different times are repeated before or after other time points. The repetitive nature of process fragments is very useful for us to analyse process changes, judge working conditions, and diagnose abnormal working conditions. Data fragments at different periods or at different time points are also like colourful clouds floating in the data space. The advancement of data over time is an important feature of industrial big data and an important aspect of big data.

Based on the Cloud diagram description and 3M characteristics of big data clusters, this section builds a set of measurement index and calculation methods to quantitatively measure the size of big data. For big data composed of multiple clouds, the size of the big data clusters is equal to the sum of the size of each piece of data cloud:

$$\|S(\Omega, T, N)\| = \sum \|S_i(\Omega_i, T_i, N_i)\| \qquad (2)$$

So, if the number of clouds is large enough, or at least one cloud is large enough, the industrial data clusters can be called as big data.

In order to describe the size of a single cloud, we can reasonably assume that each cloud $S_i(\Omega_i, T_i, N_i)$ has a compact cube denoted by $C_i$, and the cloud can be embedded in the cube compactly. The volume of a cube $C_i$ is equal to the norm of three sides.

$$\|C_i\| = \|\Omega_i\| \times \|T_i\| \times N_i \qquad (3)$$

Equation (3) is computable and easy to calculate. As long as the formula (3) is used to calculate the size of each cloud, then the formulae (2) can be used to estimate the size of the "big data". In this way, we can give an intuitive measurement of the size of the "big data".

## 3. SEVERAL TYPICAL PARADIGMS FOR INDUSTRIAL BIG DATA CLUSTER

In the industrial fields, we should not only be able to reasonably estimate the size of the data cloud, but also grasp the approaches of solving technical issues with the data cloud. Generally, the big data clusters are often widely used in three different kinds of occasions: fusion calculation so as to improve accuracy for calculations and to use all usable information for statistics inference; model correction so as to provide a more basis for prediction and support decision making; information compensation so as to bridge over gaps between fragment information. The role of big data cluster is different in some different application fields[7-9]. Correspondingly, the paradigm is also different.

### 3.1. Paradigm of Fusion Calculation

In industrial production, there are quite a large number of data information forming a piece of data cloud floating over the factory, such as the changing production process, raw material ratio data, production environment monitoring data, working condition data, and the data collected by various sensors continuously. These multi-source heterogeneous data form various types of data clouds. If these data clouds overlap, the overlapping data clouds can give us valuable measurement information of objects from different perspectives. Making full use of overlapping data clouds is helpful for us to improve the accuracy of calculation results.

Data fusion is one of the important ways of big data application. Intuitive understanding is that different data can bring different information. Combing together these data, quite a lot of information can be integrated together, which is helpful to eliminate errors and to correct prejudices. In this way, we get more and more accurate results and approximate correct inferences. In other words, data fusion technology is an information processing technology which is used to analyse various observations under certain criteria so as to complete the required

decision- making and evaluation tasks. Data fusion has achieved amazing development in the past ten years and has entered many different application areas.

In a big data environment, not every floating data cloud can participate in fusion. Data fusion is a purposeful activity. The data cloud that can participate in the fusion must be consistent in time, space, or object connotation. At the very least, if a data cloud can participate in fusion, it must be able to overlap after time traversal or space translation.

There are quite a lot of approaches for big data fusion. For example, for data layer fusion, the more widely used methods include least squares adjustment, etc.; for information layer and decision layer fusion, if all sheets of the floating data clouds $\{S_i\}$ are independent, we may use the Bayesian inference and stochastic decision making:

$$P(B|S) = \frac{\sum_i P(BS_i)}{\sum_i P(S_i)} = \sum \frac{P(S_i)}{\sum_i P(S_i)} P(B|S_i) \qquad (4)$$

where B is the event to be inferred.

### 3.2. Paradigm of Model Correction

Because of its large size, big data makes it difficult for conventional statistical methods and conventional computer data processing software to work. The combination of machine learning, big data and artificial intelligence is a measure to solve the difficult problem of big data application.

This section describes a data modelling logic that combines big data with transfer learning shown in Fig. 2:

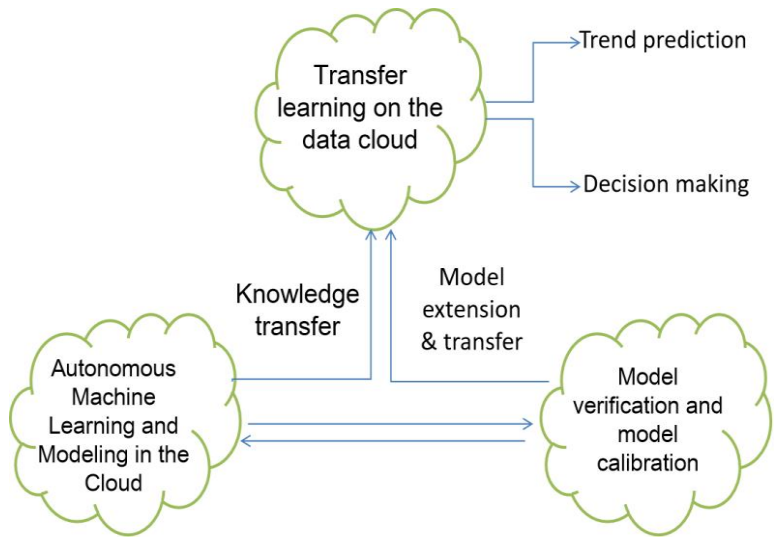

Figure 2. Prediction and Decision Based on Big Data Transfer Learning

The advantage of the above paradigm based on big data transfer learning is that the learning process is completed in the cloud chip. Due to the relative consistency of the data structure in the cloud chip, the learning process is relatively simple and easy to implement; the data of other cloud chips is used to verify the learning model And perfection to ensure the inheritance and usability of knowledge and models.

Applying big data analysis technology to specific industrial big data, the above paradigm can ensure that the big data application process is understandable and interpretable

### 3.3. Paradigm of Information Compensation

Big data is a cluster of large amounts of data. Big data has a wide range of sources, that is, multiple sources. The biggest advantage is that it can give us the opportunity to understand objects from different perspectives.

Regular-scale data gives us insight into several sides of things or objects, rather than overall impressions, similar to six blind people touching an elephant. Everyone just touches the shape of one side of the elephant. Multi-source heterogeneous data big data can be used to fill the information gaps of the untouched parts.

Based on this consideration, the fourth typical application paradigm of big data is the filling and repair of information gaps. Specifically, a map $\Psi$ is suitably constructed from the data cloud to the feature subspace:

$$\Psi: S \to \Theta \qquad (4)$$

This map has the following properties: if the big data cloud sheet $S_i$ contains the information of interest, then $\Psi(S_i) \subset \Theta$; otherwise, $\Psi(S_i) = \Phi$, an empty set. What we want to do is to find all the cloud sheets mapped to the non-empty sets from the big data cloud sheet shown in Fig.3.

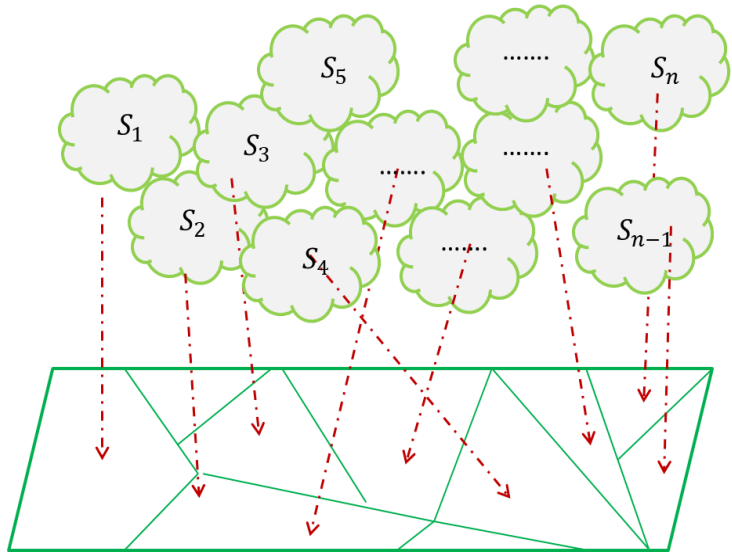

Figure 3. Mapping from big data to feature space

Obviously, if the phase space union of the big data clouds can cover the entire feature space of interest, we can use big data to achieve a complete understanding of the feature of interest. Otherwise, even if the volume of big data is large enough, it will not be possible to fully understand the changes in features from the big data cloud.

This paradigm stated above is helpful for us to adopt a problem-driven approach to delete or reduce unnecessary data blocks to achieve big data compression.

## 4. CONCLUSIONS

This paper intuitively proposes the concept of big data cloud and the 3M definition of cloud-based big data. 3M definition of big data is more suitable for industrial big data than 4V, which can fully reflect the time persistence, spatial distribution and source diversity of industrial big data. Moreover, 3M definition creates conditions for quantifying the "big" of big data. Based on 3M definition, this paper establishes a big data measure index based on the cloud compact cube volume for industrial big data, which has guiding significance and reference value for quantifying the size of big data.

On the basis of proposing big data cloud and measurement index, this paper establishes three typical application paradigms around typical applications of industrial big data, including the fusion calculation paradigm, the model correction paradigm and the information compensation paradigm, etc. These three typical application paradigms cover the basic form of the big data applications in the industrial field, which is helpful for systematically grasping the methods and approaches of industrial big data application.

ACKNOWLEDGEMENTS

This paper is financially supported by the National Nature Science Foundation of China (No.61973094, No.61933013).

## AUTHORS

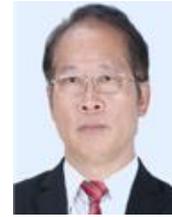

**Hu Shaolin**, professor, PhD supervisor, senior member of China Automation Society as well as China Digital Simulation Union Council; member of the editorial board of Journal of System Simulation and Journal of Ordnance and Equipment Engineering; director of the academic committee and technical director of "engineering technology center for monitoring and evaluation of smart city infrastructure" of Guangdong province. He has been committed to process monitoring, safety control, statistical learning and big data technology and application research of dynamic systems with complex structures for a long time. He has presided over 8 national natural science foundation projects and over 5 major/key pre-research projects. He has published 5 books and more than 210 papers. His recent research interests are artificial intelligence and big data technology, situational awareness and fault diagnosis, process monitoring and system safety.